\documentclass{amsart}

\usepackage{amsmath}
\usepackage{amssymb}

\theoremstyle{plain}   
\begingroup 

\newtheorem{theorem}{Theorem}[section]   
\newtheorem{corollary}[theorem]{Corollary}     
\newtheorem{lemma}[theorem]{Lemma}         
\newtheorem{proposition}[theorem]{Proposition}  
\endgroup

\theoremstyle{definition}
\newtheorem{definition}[theorem]{Definition}   

\theoremstyle{remark}
\newtheorem{remark}[theorem]{Remark}        
\newtheorem{example}[theorem]{Example}        

\newcommand{\Rn}{{\mathbb R}^{n}}

\newcommand{\R}{{\mathbb R}}
\newcommand{\N}{{\mathbb N}}

\newcommand{\Q}{\mathbb Q}

\newcommand{\lin}{\operatorname{span}}

\newcommand{\vare}{\varepsilon}
\newcommand{\vp}{\varphi}

\newcommand{\Lip}{\operatorname{Lip}}

\newcommand{\interior}{\operatorname{int}}

\newcommand{\Leb}{{\mathcal{L}}}

\newcommand{\mcA}{{\mathcal{A}}}
\newcommand{\mcC}{{\mathcal{C}}}

\begin{document}


\title{Cone monotone mappings: continuity and differentiability}

\author{Jakub Duda}

\thanks{The author was supported in part by ISF}

\email{jakub.duda@weizmann.ac.il}


\address{
Department of Mathematics, Weizmann Institute of Science,
Rehovot 76100, Israel}

\date{December 14, 2005}

\begin{abstract} 
We generalize some results of Borwein, Burke, Lewis, and Wang
to mappings with values in metric (resp.\ ordered normed linear) spaces. 
We define two classes of monotone mappings
between an ordered linear space and a metric space
(resp.\ ordered linear space): $K$-monotone dominated
and cone-to-cone monotone mappings.
$K$-monotone dominated mappings naturally generalize mappings with finite
variation (in the classical sense) and $K$-monotone functions defined
by Borwein, Burke and Lewis, to mappings with
domains and ranges of higher dimensions. First, using results
of Vesel\'y and Zaj\'\i\v{c}ek, we show some
relationships between these classes. 
Then, we show that every $K$-monotone function $f:X\to\R$,
where $X$ is any Banach space, is continuous outside of a set which can be covered
by countably many Lipschitz hypersurfaces. This sharpens a result due
to Borwein and Wang. 
As a~consequence, we obtain a~similar result for $K$-monotone dominated
and cone-to-cone monotone mappings.
Finally, we prove several
results concerning almost everywhere differentiability
(also in metric and $w^*$-senses) of these mappings.
\end{abstract}

%



\maketitle

\section{Introduction}

Let $X$ be a normed linear space, $K\subset X$ be a cone.
The cone induces ordering $\leq_K$ on $X$ as follows:
$x\leq_K y$ if $y-x\in K$.
Borwein, Burke and Lewis~\cite{BBL} 
defined the {\em $K$-increasing} functions $f:X\to\R$
as those functions that satisfy $f(x)\leq f(y)$
whenever $x\leq_K y$. They say that a function $f:X\to\R$ is
{\em $K$-monotone} provided $f$ or $-f$ is $K$-increasing.
The continuity and differentiability of $K$-monotone functions was
studied by Chabrillac, Crouzeix~\cite{CC} (in case of $\Rn$ with
the standard coordinate ordering), Borwein, Burke, Lewis~\cite{BBL},
Borwein, Wang~\cite{BW}, and others. Borwein and Goebel~\cite{BG}
provided examples showing the necessity of certain assumptions on the cone~$K$
in~\cite{BBL,BW}. The authors of~\cite{BBL} mention on page~1075 that
the Rademacher's theorem holds when the range of the function is a
Banach space with RNP, but that ``it is not clear what is true for
cone-monotone operators''. We address this issue in the current paper.
\par
We define two classes of monotone mappings with values in
metric (resp.\ ordered normed linear) spaces: $K$-monotone dominated
mappings and cone-to-cone monotone mappings.
The idea of a dominating function
in the definition of $K$-monotone dominated mappings 
is similar to the idea of a control function in the definition of d.c.\ mappings
(see~\cite{VZ}). However, the analogy should not be taken
too far as simple examples show that for instance
the $K$-monotone dominated mappings do not inherit
differentiability properties from their dominating functions
(on the other hand, d.c.\ mappings do inherit certain
differentiability properties of their control functions -- see~\cite{VZ}).
Still, continuity (resp.\ pointwise-Lipschitzness) 
of the dominating function forces the dominated
mapping to be continuous (resp.\ pointwise-Lipschitz); see Lemma~\ref{ptwiselem}.
Similarly as in the case of d.c.\ mappings, there usually
does not exist a canonical dominating function for a $K$-monotone dominated
mapping $F:X\to Y$
(except, perhaps, in the case $X=\R$).
A mapping $F:\R\to Y$ is $K$-monotone dominated if and 
only if $F$ has locally bounded variation and thus
our class generalizes mappings with bounded variation.
\par
Let us describe the structure of the present paper.
Section~\ref{Preliminaries} contains basic definitions and facts.
Section~\ref{Relationships} shows some relationships between
$K$-monotone dominated and cone-to-cone monotone mappings.
This section is motivated by results of Vesel\'y and Zaj\'\i\v{c}ek~\cite{VZoncon}.
Section~\ref{Continuity} contains results about continuity 
of monotone dominated and cone-to-cone monotone mappings.
We prove that every $K$-monotone function $f:X\to\R$,
where $X$ is an arbitrary Banach space and $K$ is a
convex cone with non-empty interior, is continuous outside of a set
which can be covered by countably many Lipschitz hypersurfaces.
This sharpens~\cite[Proposition~6]{BW} of Borwein and Wang and seems
to be new even in the case $X=\Rn$.
Then we obtain similar results for $K$-monotone dominated
and cone-to-cone monotone mappings as a consequence.
Section~\ref{Differentiability} contains
results about a.e.\ G\^ateaux differentiability of monotone dominated
and cone-to-cone monotone maps,
and section~\ref{Metrdif} contains results about metric and $w^*$-differentiability
of $K$-monotone dominated maps with values in metric spaces.
\par
Most results are formulated for mappings defined on the whole space,
however all of them can be localized; i.e.\ they hold also
for maps defined only on open sets. 
We leave the details to the interested reader.
Also, a simple example based on Lemma~\ref{lipex} and
the well known nowhere Fr\'echet differentiable
map $f:\ell_2\to\ell_2$, $f((x_n)_n)=(|x_n|)_n$ shows that there is
no hope to establish Fr\'echet differentiability of $K$-monotone
mappings between infinite-dimensional spaces. 
Every Lipschitz $f:X\to Y$ (such that $X,Y$ are Banach spaces)
is also $K$-monotone dominated (for some convex cone $K\subset X$ with non-empty interior)
by Lemma~\ref{lipex}, and thus the Lipschitz theory shows
that we can only expect G\^ateaux differentiability a.e.\ (in the appropriate
sense) of $K$-monotone dominated mappings between Banach spaces 
provided $X$ is separable, and $Y$ has RNP; see e.g.\ \cite[Chapter~6]{BL}.
Since all $K$-monotone functions are $K$-monotone dominated (the dominating function is
the function itself), we cannot expect that we will be able
to prove any results provided the cone $K$ is too small; see~\cite[Section~6]{BBL}
and~\cite{BG}. Example~\ref{mujexample} shows that if the cone in the
target space is not properly positioned, then there might be cone-to-cone monotone mappings
which are nowhere differentiable.
It remains open whether in Theorem~\ref{domdifthm} and Corollaries~\ref{domcor1p}
and~\ref{domcor2} we can replace the family $\tilde\mcC$ by $\tilde\mcA$.

\section{Preliminaries}\label{Preliminaries}

All normed (and Banach) spaces are real.
Let $X$ be a normed linear space. By $B(x,r)=B_X(x,r)$ we denote
the closed ball of $X$ (with center $x$ and radius $r$)
and by $S(x,r)=S_X(x,r)$ the corresponding sphere (omitting the subscript where no
confusion is possible). We will write $S_X=S_X(0,1)$.
\par 
We say that $X$ is an {\em ordered normed linear space} provided
it is a normed linear space equipped with an (antisymmetric) partial
ordering $\leq$ such that for $x,y,z\in X$ and $\lambda\geq0$, the 
implications
$x\leq y\implies x+z\leq y+z$,
and
$x\leq y\implies \lambda x\leq \lambda y$,
hold.
Then the corresponding cone $X_+:=\{x\in X:x\geq0\}$ is called the {\em positive cone
of $X$}.
We say that {\em $M\subset X$ has 
an upper bound} provided there exists $e\in X_+$ such that $m\leq e$
for each $m\in M$.
It is easy to see that if $X_+$ is convex,
then we have $X=X_+-X_+$ if and only if
every pair of $x,y\in X$ has an upper bound in $X_+$.
\par
We say that an ordered normed linear space $Y$ is a {\em Banach lattice}
if it is a Banach space, each pair of elements
of $Y$ has a supremum and an infimum, and $0\leq |x|\leq |y|$
implies $\|x\|\leq \|y\|$ (where $|x|:=\sup(x,-x)$).
Let $X$ be a Banach lattice. We say that $X$ has
the {\em $\sigma$-Levi property} if each norm-bounded non-decreasing
sequence in $X_+$ has the least upper bound (see~\cite{VZcit2}).
\par
Suppose that $X$ is an ordered normed linear space. 
We say that a convex subset $B\subset X$ is a {\em base}
for the cone $X_+$ if for each $y\in K\setminus\{0\}$ there exists
a unique $\lambda>0$ such that $\lambda y\in B$. Following~\cite{Jameson}
(p.\ 120), we say that $X_+$ is {\em well-based} if it has a bounded
base $B$ such that $0\not\in\overline{B}$. By~\cite[3.8.12]{Jameson},
$X_+$ is well-based if and only if
there exists $\vp\in X^*$ such that
$\vp(u)\geq \|u\|$ for each $u\in X_+$. 
\par
If $(Y,\rho)$ is a metric space, and $F:[a,b]\to Y$, then we define
the {\em variation $\bigvee^b_a f$} as a supremum of
the sums
\[ \sum_{i=1}^n \rho(f(x_i),f(x_{i-1})),\]
taken over all partitions $\{a=x_0<\dots<x_n=b\}$ of $[a,b]$.
We say that $F:\R\to Y$ has {\em locally finite variation}
provided $\bigvee^b_a f<\infty$ for all $-\infty<a<b<\infty$.
\par
\begin{definition}\label{Kmondomdef}
Let $X$ be a normed linear
space, $K\subset X$ be a non-empty cone, let $(Y,\rho)$ be a metric space. 
We say that $F:X\to Y$ is {\em $K$-monotone dominated} provided
there exists $h:X\to\R$ such that 
\begin{equation*}\label{mondomdef} 
\rho(H(x+k),H(x))\leq h(x+k)-h(x),
\end{equation*}
whenever $x\in X$ and $k\in K$. 
Then we say that $h$ is {\em the dominating function for $H$}.
\end{definition}
Trivially,
every dominating function is $K$-monotone.
If we fix $K,X,Y$, and assume that $Y$ is a normed linear space, 
then it is easy to see that $K$-monotone mappings $F:X\to Y$ form a vector space.
Since every metric space embeds isometrically into some $\ell_\infty(\Gamma)$
(for some $\Gamma$),
there would be no loss of generality in assuming that the space $Y$
in the definition of $K$-monotone dominated mappings is a normed linear space.
\begin{lemma}\label{lipex}
Let $X$ be a normed linear space, $(Y,d)$ be a metric space.
If $f:X\to Y$ is Lipschitz, then $f$ is $K$-monotone dominated for
some convex cone $K\subset X$ with non-empty interior.
\end{lemma}
\begin{proof} Take $0\neq x\in X$ with $\|x\|=1$. Find $x^*\in X^*$ such that
$\|x^*\|=\langle x^*,x\rangle=1$. Let $\alpha\in(0,1)$, and put $K:=\{y\in X:\alpha\|y\|\leq x^*(y)\}$.
Then $K$ is a convex cone with non-empty interior. If $y,z\in X$ are such that $y-z\in K$,
then $d(f(y),f(z))\leq L\|y-z\|\leq\frac{L}{\alpha}\cdot x^*(y-z)$, and thus $\frac{L}{\alpha}\cdot x^*$ is a dominating
function for $f$, and $f$ is $K$-monotone dominated.
\end{proof}
If $X=\R$, then $\R_+$-monotone dominated mappings are
exactly those mappings that have locally finite variation; see Section~\ref{Relationships}
for more details.
\par
\begin{definition}\label{conetoconedef}
Let $X$, $Y$ be ordered normed linear spaces, $F:X\to Y$. 
We say that $F$ is {\em $(X_+,Y_+)$-increasing} provided 
$x\leq y\implies F(x)\leq F(y)$,
whenever $x,y\in X$.
We say that $F$ is {\em $(X_+,Y_+)$-decreasing} provided $-F$ is
$(X_+,Y_+)$-increasing. Finally, we say that $F$ is {\em $(X_+,Y_+)$-monotone}
provided it is either $(X_+,Y_+)$-increasing or  it is $(X_+,Y_+)$-decreasing.
\end{definition}
We shall also refer to $(X_+,Y_+)$-monotone mappings as
cone-to-cone monotone mappings.
\par
We say that $F:X\to Y$ is 
{\em pointwise-Lipschitz at $x\in X$} provided
$\Lip(F,x):=\limsup_{t\to x}\frac{\|F(t)-F(x)\|}{\|t-x\|}<\infty$.
We say that $F:X\to Y$ is {\em Lipschitz} provided
there exists $C>0$ such that $\|F(x)-F(y)\|\leq C\|x-y\|$
whenever $x,y\in X$.
We say that $F:X\to Y$ is {\em G\^ateaux differentiable at $x\in X$}
provided
$T(h):=\lim_{t\to0} \frac{F(x+th)-F(x)}{t}$
exists for each $h\in X$, and $T:X\to Y$ is a bounded linear
operator. For basic facts about G\^ateaux derivatives,
see~\cite{BL}.
Given $f:X\to\R$, we use $\overline{f}$ (resp.\ $\underline{f}$)
for the {\em upper (resp.\ lower) semi-continuous envelope} of $f$.
\par
Let $X,Y$ be (real) normed linear spaces.
For $f:X\to Y$ we shall denote 
\begin{equation*}
MD(f,x)(u)=\lim_{r\to0}\frac{\|f(x+ru)-f(x)\|}{|r|}\quad\text{ for }x,u\in X.
\end{equation*}
It was defined in~\cite{K}.
If $MD(f,x)(u)$ exists for all $u\in X$, we say that
{\em $f$ is directionally metrically differentiable at $x$}. 
We will say that {\em $f$ is metrically G\^ateaux differentiable at $x$}
provided $f$ is directionally metrically differentiable at $x$, and
$MD(f,x)(\cdot)$ is a continuous seminorm.
We say that {\em $f$ is metrically differentiable at $x$}, 
provided $f$ is metrically G\^ateaux differentiable, and 
\begin{equation}\label{mdf}
\|f(z)-f(y)\|-MD(f,x)(z-y)=o(\|z-x\|+\|y-x\|),\quad\text{ when }(y,z)\to (x,x).
\end{equation}
\par
We will also need the notion $w^*$-G\^ateaux derivatives. 
It goes back to~\cite{HM}.
Let $X,Y$ be separable Banach spaces, $f:X\to Y^*$ be a mapping. For $v\in X$ we say
that {\em $wd(f,x)(v)$ exists } provided
$wd(f,x)(v)=w^*-\lim_{t\to0}\frac{f(x+tv)-f(x)}{t}$
exists. 
We say that {\em $f$ is $w^*$-G\^ateaux differentiable at $x$} provided
$wd(f,x)(v)$ exists for all $v\in X$, and $wd(f,x)(\cdot)$ is a bounded linear map.
We say that {\em $f$ is $w^*$-Fr\'echet differentiable at $x$} provided
$f$ is $w^*$-G\^ateaux differentiable at $x$, and 
\begin{equation}\label{wstardif} 
w^*-\lim_{y\to x} \frac{f(y)-f(x)-wd(f,x)(y-x)}{\|y-x\|}=0.
\end{equation}
$w^*$-G\^ateaux (and $w^*$-Fr\'echet) differentiability of pointwise
Lipschitz mappings was studied in~\cite{Dmwd}, where we introduced these notions.
\par
Suppose that $X$ is a normed linear space. We say that 
$M\subset X$ is a {\em Lipschitz hypersurface} provided
there exists $x^*\in X^*$, a Lipschitz function $f:Y=\{x^*=0\}\to\R$,
and $0\neq v\subset X$
such that $M=\{y+f(y)v:y\in Y\}$.
We say that $A\subset X$ {\em can be covered by countably
many Lipschitz hypersurfaces} provided there exist Lipschitz
hypersurfaces $L_n\subset X$ such that $A\subset \bigcup_n L_n$.
\par
Preiss and Zaj\'\i\v{c}ek
introduced the $\sigma$-directionally porous sets in~\cite{PZ}.
Note that every Lipschitz hypersurface is a directionally porous
set, and thus every set which can be covered by countably
many Lipschitz hypersurfaces is $\sigma$-directionally porous.
For a recent survey of use of negligible sets in similar
contexts, see~\cite{ZaA}.
\par
Let $X$ be a separable Banach space,
$A\subset X$, and $0\neq u\in X$. 
We say that {\em $A\in{\mathcal{A}}(u)$} provided $A$ is Borel and
${\mathcal{L}}^1(\{\lambda\in\R:x+\lambda u\in A\})=0$,
for all $x\in X$.
For a sequence $\{u_n\}\subset X$ we define
\[{\mathcal{A}}(\{u_n\})=\bigg\{E\in X:E=\bigcup_n E_n\text{ with }E_n\in{\mathcal{A}}(u_n)\bigg\}.\]
Finally, we say that {\em $A$ is Aronszajn null} provided $A$ is Borel
and for each complete sequence $\{u_n\}\subset X$ we have $A\in{\mathcal{A}}(\{u_n\})$
(a sequence $\{u_n\}$ is {\em complete} provided $X={\overline{\lin}}(\{u_n\})$).
For more information about Aronszajn null sets, see~\cite{BL}.
\par
Let $X$ be a separable Banach space.
If $0\neq v\in X$, 
then let $\tilde\mcA(v,\vare)$ be the system of all Borel sets $B\subset X$
such that $\{t:\vp(t)\in B\}$ is Lebesgue null whenever $\vp:\R\to X$ is such
that the function $t\to\vp(t)-tv$ has Lipschitz constant at most~$\vare$, 
and $\tilde\mcA(v)$ is the system of all sets $B$ such that
$B=\bigcup_{k=1}^\infty B_k$, where $B_k\in\tilde\mcA(v,\vare_k)$ for some $\vare_k>0$.
\par
We define $\tilde\mcA$ (resp.\ $\tilde\mcC$) as the system of those $B\subset X$ that can be,
for all complete sequences $\{v_n\}$ in $X$ (resp.\
for some sequence $\{v_n\}$ in $X$), written
as $B=\bigcup_{n=1}^{\infty} B_n$, where each $B_n$ belongs
to $\tilde\mcA(v_n)$.
These families of sets were defined in~\cite{PZ} as proper subfamilies
of Aronszajn null sets.
\par
We will also use the following notation: The symbol $C_b(T)$
(resp.\ $C(K)$) denotes the space of all continuous bounded
functions on an arbitrary topological space $T$ (resp.\ of
all continuous functions on an arbitrary compact space $K$)
equipped with the supremum norm. When we deal with the spaces
$L^p(\mu)$, we allow an arbitrary measure $\mu$.
\par
The following lemma shows that the cone monotone mappings
do inherit some properties of their dominating functions
(see also Introduction).

\begin{lemma}\label{ptwiselem}
Let $X,Y$ be a normed linear spaces, let $K\subset X$ be a convex cone with
non-empty interior.
Suppose that $H:X\to  Y$ is $K$-monotone dominated with a~dominating
function $h:X\to\R$. 
\begin{enumerate}
\item
If $h$ is continuous at $x\in X$,
then $H$ is continuous at $x$.
\item
If $h$ is pointwise-Lipschitz at $x\in X$,
then $H$ is pointwise Lipschitz at $x$.
\end{enumerate}
\end{lemma}

\begin{proof}
We will only prove~(ii), as the proof of~(i) is similar.
Without any loss of generality, we can assume that $x=0$.
Because $K$ has a non-empty interior, choose $k\in K$ such that
\begin{equation}\label{choicek} 
k+y\in K\quad\text{ whenever }y\in B(0,1).
\end{equation}
Put $C:=2\|k\|$. Then it is easy to see that each $y\in X$
can be written as $y=k_1-k_2$, where $k_1,k_2\in K\cap B(0,C\|y\|)$.
For any $\vare>0$ find $\delta>0$ such that 
$|h(t)-h(0)|\leq (1+\vare)\Lip(h,0)|t|$
for any $t\in B(0,\delta)$.
Suppose that $y\in B(0,\delta/C)$. Then
there exist $k_1,k_2\in K\cap B(0,\delta)$ such that $y=k_1-k_2$,
and $\max(\|k_1\|,\|k_2\|)\leq C\|y\|$.
Now,
\begin{equation}\label{j1}
\begin{split} 
\|H(y)-H(0)\|&=\|H(k_1-k_2)-H(0)\|\\
&\leq\|H(k_1-k_2)-H(-k_2)\|+\|H(-k_2)-H(0)\|\\
&\leq h(k_1-k_2)-h(-k_2)+h(0)-h(-k_2),
\end{split}
\end{equation}
and
\begin{equation}\label{j2} 
\begin{split}
\|H(y)-H(0)\|&=\|H(k_1-k_2)-H(0)\|\\
&\leq\|H(k_1-k_2)-H(k_1)\|+\|H(k_1)-H(0)\|\\
&\leq h(k_1)-h(k_1-k_2)+h(k_1)-h(0).
\end{split}
\end{equation}
By adding~\eqref{j1} and~\eqref{j2}, we obtain
\begin{equation*}
\begin{split}
 2\|H(y)-H(0)\|&\leq 2 (h(k_1)-h(-k_2))\\
&\leq 2(h(k_1)-h(0)+h(0)-h(-k_2))\\
&\leq 2(1+\vare)\,C\,\Lip(h,0)\|y\|.
\end{split}
\end{equation*}
By sending $\vare\to0$ we have $\Lip(H,x)\leq C\Lip(h,x)$, and the conclusion follows.
\end{proof}

The following lemma shows that for
$K$-monotone functions G\^ateaux differentiability
implies pointwise-Lipschitzness.

\begin{lemma}\label{pttoptlem}
Let $X$ be a normed linear space, $K\subset X$ a closed
convex cone with non-empty interior, and $f:X\to\R$ be $K$-monotone.
If $f$ is G\^ateaux differentiable at $x$,
then $f$ is pointwise-Lipschitz at $x$.
\end{lemma}

\begin{proof}
Let $k\in K$ be as in~\eqref{choicek}.
Let $C>0$ be such that each $y\in X$ can be written
as $y=k_1-k_2$ with $k_1,k_2\in K\cap B(0,C\|y\|)$
(see the text after~\eqref{choicek}).
Without any loss of generality, we can assume that $x=0$. 
Let $\delta>0$ be such that $|f(\lambda k)-f(0)|\leq M|\lambda|$
whenever $|\lambda|<\delta$.
Let $y\in B(0, \delta/(C\|k\|))$. There exist $k_1, k_2\in K\cap B(0,\delta/\|k\|)$ such that 
$y=k_1-k_2$, and $\max(\|k_1\|,\|k_2\|)\leq C\|y\|$.
Since 
\begin{equation}\label{f1} 
|f(y)-f(0)|\leq f(k_1)-f(k_1-k_2)+f(k_1)-f(0),
\end{equation}
and
\begin{equation}\label{f2} 
|f(y)-f(0)|\leq f(k_1-k_2)-f(-k_2)+f(0)-f(-k_2),
\end{equation}
by adding~\eqref{f1} to~\eqref{f2}, we obtain
\begin{equation}\label{f3} |f(y)-f(0)|\leq f(k_1)-f(-k_2).\end{equation}
By~\eqref{choicek}, we have that $-k_i+\|k_i\|k\in K$ for $i=1,2$. 
Thus
\begin{equation}\label{f4} 
f(k_1)-f(0)\leq f(\|k_1\|k)-f(0)\leq f(C\|y\|k)-f(0)
\leq CM \|y\|.
\end{equation}
Similarly,
\begin{equation}\label{f5} 
-f(-k_2)+f(0)\leq -f(-\|k_2\|k)+f(0)\leq -f(-C\|y\|k)+f(0)\leq CM \|y\|.
\end{equation}
Putting~\eqref{f3},~\eqref{f4}, and~\eqref{f5} together, we obtain that 
$f$ is pointwise Lipschitz at $x$.
\end{proof}

The following auxiliary proposition tells us that
the composition of a cone-to-cone monotone mapping with
a $K$-monotone dominated mapping is again $K$-monotone dominated
provided the intermediate cone is the same.

\begin{proposition}\label{comp-prop}
Let $X$, $Y$ be ordered normed linear spaces, let $Z$ be a linear space.
Suppose that $G:X\to Y$ is $(X_+,Y_+)$-monotone, and
that $H:Y\to Z$ is $Y_+$-monotone dominated.
Then $H\circ G$ is $X_+$-monotone dominated.
\end{proposition}

\begin{proof} 
Without any loss of generality, we
can assume that $G$ is $(X_+,Y_+)$-increasing.
Suppose that $x\leq y$ for some $x,y\in X$.
Then $G(x)\leq G(y)$, and thus
$h(G(x))\leq h(G(y))$. We have that $H\circ G$ is $X_+$-monotone
dominated with the dominating function $h\circ G$.
\end{proof}

\section{Relationships between the classes of cone monotone mappings}
\label{Relationships}

%

If $Y$ is a normed linear space, and $F:\R\to Y$,
then it is easy to see that
$F$ is $\R_+$-monotone dominated mapping if and only
if $F$ has locally finite variation. These mappings were 
studied before; see e.g.\ \cite{VZoncon} for references.
\par
The following theorem is proved in~\cite[Theorem~2.7]{VZoncon}:

\begin{theorem}\label{vzthm}
Let $I\subset\R$ be an open (or closed) interval, $Y$ be a Banach lattice
with the $\sigma$-Levi property, and $f:I\to Y$ be a mapping having
locally finite variation.
\par
Then there exist nondecreasing mappings $g,h:I\to Y$ such that
$f=g-h$ and $g,h$ have locally finite variation. Moreover,
the decomposition $f=g-h$ is minimal in the class of all representations
of $f$ as the difference of nondecreasing mappings, i.e.: if $f=g^*-h^*$
is such a representation then $g(\beta)-g(\alpha)\leq g^*(\beta)-g^*(\alpha)$
for all $\alpha<\beta$, $\alpha,\beta\in I$.
\end{theorem}

It has the following corollary:

\begin{corollary}\label{coro10}
Let $Y$ be a Banach lattice with the $\sigma$-Levi property,
and $f:\R\to Y$ be $\R_+$-monotone dominated. Then
$f$ can be written as a difference of two $(\R_+,Y_+)$-monotone
mappings.
\end{corollary}

The following is a direct consequence of Corollary~\ref{coro10}
and~\cite[Remark~2.8]{VZoncon}; we formulate it for the reader's
convenience:

\begin{corollary}\label{coro11}
Let $Y$ be any dual Banach lattice (in particular
$L_p(\mu)$, $1<p<\infty$, $\ell_\infty(\Gamma)$ for any $\Gamma$,
or $L_\infty$ when $\mu$ is $\sigma$-finite).
If $F:\R\to Y$ is $\R_+$-monotone dominated, then
$F$ can be written as a difference of two $(\R_+,Y_+)$-monotone
mappings.
\end{corollary}

The following proposition gives a general condition when
$K$-monotone dominated mappings between two ordered normed linear
spaces can be written as differences of cone-to-cone monotone mappings.

\begin{proposition}\label{unitprop}
Let $Y$ be an ordered normed linear space whose
unit ball has an upper bound, let $X$ be an ordered normed
linear space. Then each for each $X_+$-monotone dominated mapping
$F:X\to Y$ can be written as a difference of two $(X_+,Y_+)$-monotone mappings.
\end{proposition}

\begin{proof}
Let $e\in Y_+$ be the upper bound of $B_Y$, i.e.\ $y\leq\|y\|e$
for each $y\in Y$, and let $f:X\to\R$ be the dominating function for $F$.
Then for $x,y\in X$ with $x\leq y$ we have
$F(x)-F(y)\leq \|F(y)-F(x)\|e\leq f(y)e-f(x)e$.
If we put $H=F+f\cdot e$, then $H$ is $(X_+,Y_+)$-monotone by
the above, and $F=H-f\cdot e$. This gives us the conclusion
since $f\cdot e$ is trivially $(X_+,Y_+)$-monotone.
\end{proof}

\begin{corollary}\label{vzcor1}
Let $(Y,K)$ be any of the spaces $C_b(T)$, $L^\infty(\mu)$,
$C_b(T)^{**}$ or $L^\infty(\mu)^{**}$ with its canonical cone.
Let $X$ be an ordered normed linear space. 
Then every $X_+$-monotone dominated mapping $F:X\to Y$ can
be written as a difference of two $(X_+,Y_+)$-monotone mappings.
\end{corollary}

Now we will show that in some cases cone-to-cone
monotone mappings are, in fact, $K$-monotone dominated.

\begin{proposition}\label{orderimpliesmon}
Let $Y$ be an ordered normed linear space
whose cone $Y_+$ is well-based, let $X$ be
an ordered linear space. Then every $(X_+,Y_+)$-monotone
mapping $F:X\to Y$ is $X_+$-monotone dominated.
\end{proposition}

\begin{proof}
Without any loss of generality, we can assume that $F$ is 
$(X_+,Y_+)$-increasing.
Since $Y_+$ is well-based, there exists $\vp\in Y^*$ such that
$\vp(u)\geq \|u\|$ for each $u\in Y_+$. 
This shows that $Id_Y:Y\to Y$ is $Y_+$-monotone dominated.
Thus Proposition~\ref{comp-prop} implies that
$F=Id_Y\circ F$ is $X_+$-monotone dominated.
\end{proof}

\begin{corollary}\label{vzcor2}
Let $Y$ be any of the spaces $L_1(\mu)$, $L^{\infty}(\mu)^*$, $C_b(T)^*$
with its canonical cone. Let $X$ be an ordered normed linear space.
Then each $(X_+,Y_+)$-monotone mapping is $X_+$-monotone dominated.
\end{corollary}

\begin{proof} 
In the proof of~\cite[Corollary~3.2]{VZoncon} it
is shown that the canonical cone of each of the spaces is
well-based. Thus Proposition~\ref{orderimpliesmon} gives the conclusion.
\end{proof}

\section{Continuity of cone monotone mappings}\label{Continuity}

Let $(X,\|\cdot\|)$ be a normed linear space. We say that $\|\cdot\|$ is {\em LUR at $x\in S_X$}
provided $x_n\to x$ whenever $\|x_n\|=1$, and $\|x_n+x\|\to 2$.
We say that $\|\cdot\|$ is {\em LUR} (or {\em locally
uniformly rotund}) provided it is LUR at each point $x\in S_X$.
For more information about rotundity and renormings, see~\cite{DGZ}.
\par
We will need the following renorming result, which is 
proved e.g.\ in~\cite[Lemma II.8.1]{DGZ}.

\begin{lemma}\label{renormlem} 
Let $Y$ be a subspace of a Banach space $X$
and let $|\cdot|$ be an equivalent norm on~$Y$. Then $|\cdot|$
can be extended to an equivalent norm on~$X$. If $|\cdot|$
is an equivalent locally uniformly rotund norm on~$Y$ then
$|\cdot|$ can be extended to an equivalent norm on $X$ which
is locally uniformly rotund at each point of $Y$.
\end{lemma}

\begin{remark}\label{renormrem} 
Let $(X,\|\cdot\|_X)$ be a normed linear space, $y\in S_X$.
By $\hat X$ denote the Banach space completion of $X$
(i.e.\ $\hat X=\overline X$, $\|x\|_{\hat X}=\|x\|_X$ 
for each $x\in X$, and $\hat X$ is a Banach space).
Then we can apply Lemma~\ref{renormlem} to $\hat X$ and $Y=\lin\{y\}$
to obtain an equivalent norm $\|\cdot\|_1$ on $\hat X$ which 
is LUR at each point of $Y$, and with the property that $\|\lambda y\|_1=\|\lambda y\|_X$
for each $\lambda\in\R$. 
It follows that $\|\cdot\|_1$ (restricted to $X$) is LUR at each point of $Y$ on $X$.
Thus Lemma~\ref{renormlem} also holds
for any normed linear space $X$ and $Y=\lin\{y\}$, where $y\in S_X$.
\end{remark}

Using similar reasoning as in Zaj\'\i\v{c}ek~\cite[Lemma~1]{Zmult}
together with a renorming Lemma~\ref{renormlem} (resp.\ Remark~\ref{renormrem}),
we can improve~\cite[Proposition~12]{BW}, 
since it is easy to see that
every Lipschitz hypersurface is a directionally porous set.

\begin{proposition}\label{contprop}
Let $X$ be a normed linear space. Assume that $K\subset X$ 
is a convex cone with non-empty interior and $f:X\to\R$ is $K$-monotone.
Then $D:=\{x\in X:f\text{ is discontinuous at }x\}$
can be covered by countably many Lipschitz hypersurfaces.
\end{proposition}

\begin{proof}
Following the proof of~\cite[Proposition~12]{BW}
we have $D=\{x\in X:\underline{f}(x)<\overline{f}(x)\}$.
Write 
$S_1:=\{x\in X:\underline{f}(x)<f(x)\}$, and 
$S_2:=\{x\in X:{f}(x)<\overline{f}(x)\}$.
We will only prove that $S_2$ can be covered by countably
many Lipschitz hypersurfaces (the proof for $S_1$ is similar).
Write $S_2=\bigcup_{p\in\Q} D_p$, where
$D_p:=\{x\in X:f(x)<p<\overline{f}(x)\}$.
As in the proof of~\cite[Proposition~12]{BW}, we see that
\begin{equation}\label{asin} 
[x-\interior(K)]\cap D_p=\emptyset,
\end{equation}
for each $x\in D_p$. 
\par
We will prove that each $D_p$ is contained in a Lipschitz hypersurface.
Fix $p\in\Q$.
Choose $k\in\interior(K)$. We can assume that $\|k\|=1$.
By Remark~\ref{renormrem}, find an equivalent norm $\|\cdot\|_1$ on $X$
such that $\|\cdot\|_1$ is LUR at~$k$ (and with $\|k\|_1=1$).
From now on, we consider $X$ equipped with $\|\cdot\|_1$.
Then $B(k,\vare)\subset K$ for some $\vare>0$.
Choose $x^*\in X^*=(X,\|\cdot\|_1)^*$ such that $\|x^*\|_1^*=\langle x^*,k\rangle=1$.
Since $\|\cdot\|_1$ is LUR at $k$,
we have  that if $x_n\in X$ are such that $\|k+x_n\|_1\to 2$,
and $\|x_n\|_1=1$, then $x_n\to k$.
Let $0<\alpha<1$. Then $K_\alpha=\{x\in X:\alpha\|x\|_1< x^*(x)\}$
is a convex cone with non-empty interior.
If $x\in K_\alpha\cap S_X$, then
$1+\alpha=x^*(k+x)\leq \|k+x\|_1$,
and by the LUR property of $\|\cdot\|_1$ at~$k$,
we see that there exists $0<\alpha<1$ such that
$K_\alpha\cap S_X\subset B(k,\vare/2)\subset\interior(K)$ (and thus $K_\alpha\subset\interior(K)$).
\par
By~\eqref{asin},
we have
$[x-\interior(K)]\cap D_p=\emptyset$,
for any $x\in D_p$, and thus
$[x-K_\alpha]\cap D_p=\emptyset$,
for any $x\in D_p$.
Let $Y:=\{x^*=0\}$.
We have that if $x,y\in D_p$, then 
$ |x^*(x-y)|\leq\alpha\|x-y\|_1$.
This implies that if $x,y\in D_p$, then
\begin{equation}\label{pieq} 
\begin{split}
(1-\alpha)|x^*(y)-x^*(x)|&\leq(1-\alpha)\|y-x\|_1\\
&\leq \|x-y\|_1-\alpha|x^*(y-x)|\\
&\leq \|x-y\|_1-\|(x^*(y-x))k\|_1\\
&\leq \|x-y-x^*(x-y)k\|_1.
\end{split}
\end{equation}
It is easy to see that the mapping $\pi:D_p\to Y$
defined as $\pi(x):=x-x^*(x)k$ is one-to-one.
Suppose that $x',y'\in \pi(D_p)$.
Then there exist $x,y\in D_p$ such that $\pi(x)=x'$,
and $\pi(y)=y'$. The inequality~\eqref{pieq} shows
that the function $f:\pi(D_p)\to\R$ defined
as $f(x')=x^*(x)$, where $x$ is
the (unique) vector in $D_p$ such that $\pi(x)=x'$,
is Lipschitz. Thus $D_p=\{y+f(y)k:y\in\pi(D_p)\}$.
Since $f:\pi(D_p)\to\R$ is Lipschitz, it can
be extended to a Lipschitz function (call it again $f$) on~$Y$.
Define $L:=\{y+f(y)k:y\in Y\}$. Then $L$ is a Lipschitz
hypersurface, and $D_p\subset L$.
\end{proof}

\begin{theorem}\label{contthm}
Suppose that $X$ is an ordered normed linear space
such that $X_+$ is convex with non-empty interior,
$(Y,d)$ is a metric space, 
and $F:X\to Y$ is $X_+$-monotone dominated.
Then the set $D$ of points where $F$ is not continuous
can be covered by countably many Lipschitz hypersurfaces.
\end{theorem}

\begin{proof}
We can isometrically embed $Y$ into $\ell_\infty(\Gamma)$ for some $\Gamma$,
and thus without any loss of generality, we can assume that $Y$ is a normed linear space.
Let $f:X\to\R$ be the dominating function for $F$.
By Proposition~\ref{contprop}, the set of points $D'\subset X$
where $f$ is discontinuous can be covered by countably many Lipschitz hypersurfaces.
Lemma~\ref{ptwiselem}(i) implies that $D\subset D'$.
\end{proof}

\begin{corollary}\label{contcor}
Suppose that $X$, $Y$ are ordered normed
linear spaces such that $X_+$ is convex 
with non-empty interior, and $Y_+$ is well-based.
If $F:X\to Y$
is $(X_+,Y_+)$-monotone, then
$D:=\{x\in X:F\text{ is discontinuous at }x\}$
can be covered by countably many Lipschitz hypersurfaces.
\end{corollary}

\begin{proof}
Proposition~\ref{orderimpliesmon} shows that $F$ is $X_+$-monotone dominated
since $Y_+$ is well-based. The conclusion now follows from Theorem~\ref{contthm}.
\end{proof}

\section{Differentiability of cone monotone mappings}\label{Differentiability}

We have the following theorem concerning G\^ateaux differentiability
of monotone dominated mappings. It is a corollary of a general version
of Stepanoff's theorem which was proved in~\cite{Dongat} 
(as a strengthening of a result due to Bongiorno~\cite{B}). 

\begin{theorem}\label{domdifthm}
Suppose that $X$ is a separable Banach space, 
$K\subset X$ is a convex cone with non-empty interior,  
$Y$ is a Banach space with RNP, and $F:X\to Y$ is $K$-monotone dominated.
Then there exists a set $A\in\tilde\mcC$ such that
$F$ is G\^ateaux differentiable at all $x\in X\setminus A$.
\end{theorem}

\begin{proof}
Let $f:X\to\R$ be a dominating function for $F$. 
Then by~\cite[Theorem~15]{Dongat},
there exists a set $A_1\in\tilde\mcC$ such that
$f$ is G\^ateaux differentiable at all $x\in X\setminus A_1$.
Lemma~\ref{pttoptlem} implies that $h$ is pointwise-Lipschitz
at each $x\in X\setminus A_1$.
By Lemma~\ref{ptwiselem} it follows that $F$ is pointwise-Lipschitz at
all~$x\in X\setminus A_1$. By~\cite[Theorem~10]{Dongat} there exists a set 
$A_2\in\tilde\mcA$ such that $F$ is G\^ateaux differentiable at all $x\in X\setminus (A_1\cup A_2)$.
Putting $A:=A_1\cup A_2$ finishes the proof of the theorem.
\end{proof}

\begin{remark} 
Using the methods of~\cite{Dongat}, we can actually obtain
a version of~\cite[Theorem~10]{Dongat} where the notion of ``G\^ateaux differentiability''
is replaced by ``Hadamard differentiability'' and thus Theorem~\ref{domdifthm}
also holds for this notion. We do not enter this subject here.
\par
If $X=\Rn$ in the previous theorem, then	we can conclude
that $F$ is even almost everywhere Fr\'echet differentiable
(the proof is analogous to the proof of Theorem~\ref{domdifthm};
it uses~\cite[Theorem~2.10]{Dmwd} and
the fact that $K$-monotone functions on $\Rn$ are (Fr\'echet) differentiable almost everywhere;
see e.g.\ \cite{BBL}).
\end{remark}
	
Using the same reasoning as in~\cite[Corollary~11]{Dongat}, we can obtain
the following corollary. It is a consequence of a recent result of Zaj\'\i\v{c}ek~\cite{Znew}
who proved that the sets in~$\tilde\mcC$ (in a separable Banach space)
are $\Gamma$-null (in the sense of~\cite{LP}) and the resutls of~\cite{LP}.

\begin{corollary}\label{domcor1p}
Let $X$ be a Banach space such that $X^*$ is separable,
$K\subset X$ be a convex cone with non-empty interior,
let $Y$ be a Banach space with RNP, $f:X\to Y$ be $K$-monotone dominated,
$g:X\to\R$ continuous convex.
Then there exists $x\in X$ such that $f$ is G\^ateaux differentiable at~$x$ 
and $g$ is Fr\'echet differentiable at~$x$.
\end{corollary}

The following corollary tells us that if the
cone in the target space is not too big, then
cone-to-cone monotone mappings are G\^ateaux differentiable
almost everywhere.

\begin{corollary}\label{domcor2}
Suppose that $X$ be an ordered separable  Banach space such that
$X_+$ is convex with non-empty interior, $Y$ be 
an~ordered Banach space with RNP such that $Y_+$ is well-based, 
and $F:X\to Y$ be $(X_+,Y_+)$-monotone.
Then there exists a set $A\in\tilde\mcC$ 
such that $F$ is G\^ateaux differentiable at all $x\in X\setminus A$.
\end{corollary}

\begin{proof}
Proposition~\ref{orderimpliesmon} implies that $F$ is $X_+$-monotone dominated.
Theorem~\ref{domdifthm} now implies that there exists a set $A\in\tilde\mcC$
such that $F$ is G\^ateaux differentiable at all $x\in X\setminus A$.
\end{proof}


The following example shows that if the cone in the target space is not
properly positioned, then Corollary~\ref{domcor2} does not hold.

\begin{example}\label{mujexample}
Let $K$ be the non-negative cone in $\ell_2$
(i.e.\ $x=(x_n)_n\in K$ iff $x_n\geq0$ for all $n\in\N$).
Then there exists a mapping $f:\R\to\ell_2$, which is $(\R_+,K)$-monotone,
but nowhere differentiable.
\end{example}

\begin{remark} In our example, the domain of $f$ is $\R$ and thus
the notions of G\^ateaux and Fr\'echet differentiability coincide.
\end{remark}

\begin{proof} 
Let $(e_n)_n$ be the canonical orthonormal basis of $\ell_2$.
We will use the following easy observation:
if $f:\R\to\ell_2$ is differentiable at~$x$, then each $f_j$
is differentiable at~$x$ (where $f(x)=\sum_j f_j(x)\cdot e_j$),
and $\max_j|f_j'(x)|\leq\|f'(x)\|$.
\par
It is easy to see that there exist  sequences $(a_n)_{n\in\N}\subset\R$
and $(m_n)_{n\in\N}\subset\N$  such that 
\begin{itemize} 
\item $a_n>0$ for all $n\in\N$, 
\item $(m_n)_n$ is increasing,
\item $\sum_n a_n^2<\infty$,
\item $\sum_{j=m_n}^{m_{n+1}-1} a_j=2n$ (for $n=1,2,\dots$).
\end{itemize}
Find $k_n>0$ such that $\sum_n (k_n\cdot a_n)^2<\infty$ and $\lim_{n\to\infty}k_n=\infty$.
Define $f_j:\R\to\R$ as 
\[f_j(x)=k_j\cdot\max\bigg(0,\min\bigg(a_j,x+n-\sum_{k=m_n}^{j-1} a_k\bigg)\bigg),\]
for $j=m_n+1,\dots,m_{n+1}-1$, and 
$f_j(x)=k_j\cdot\max(0,\min(x+n,a_j))$ for $j=m_n$.
Then we have $0\leq f_j(x)\leq k_j\cdot a_j$ for each $x\in\R$.
For $x\in\R$ define $f(x):=\sum_j f_j(x)\cdot e_j$.
It is easy to see that $f$ is well-defined, continuous, and $(\R_+,K)$-increasing
(since each $f_j$ is increasing).
Assume that $f$ is differentiable at some point $x\in\R$. 
Take any $n\in\N$ such that~$|x|<n$.
If $x=-n+\sum_{j=m_n}^k a_j$ for some
$k\in\{m_n,\dots,m_{n+1}-2\}$, then $f_k$ is
not differentiable at $x$ (and thus $f$ cannot be differentiable
at~$x$). Thus $x\in(-n,-n+a_{m_n})$ or $x\in(-n+\sum^{l}_{j=m_n}a_j,-n+\sum^{l+1}_{j=m_n}a_j)$
for some $l\in\{m_n,\dots,m_{n+1}-2\}$. But then $f'_{m_n}(x)=k_{m_n}$ or $f_l'(x)=k_l$, and since $k_j\to\infty$,
we have that $f$ is not differentiable at~$x$.
\end{proof}

\section{Metric differentiability of cone monotone mappings}
\label{Metrdif}

\begin{theorem}\label{rnmetrthm}
Suppose that $K\subset\R^n$ is a convex
cone with non-empty interior.
Let $(M,\rho)$ be a metric space, and
let $F:\R^n\to M$ be $K$-monotone dominated.
Then for almost every $x\in\Rn$ we
have that $F$ is metrically differentiable at~$x$.
\par
If $M=Y^*$, where $Y$ is a separable Banach space,
then for almost every $x\in\Rn$ we
have that $F$ is metrically differentiable at~$x$,
$F$ is $w^*$-Fr\'echet differentiable at~$x$,
and $MD(F,x)(w)=\|wd(F,x)(w)\|$ for all $w\in\Rn$.
\end{theorem}

\begin{proof}
If $M$ is a metric space, then we can embed $(M,\rho)$ isometrically into~$\ell_\infty(\Gamma)$
for some $\Gamma$, and so we can assume that $F:\Rn\to\ell_{\infty}(\Gamma)$.
Let $f:\Rn\to\R$ be a dominating function for $F$.
Then~\cite[Theorem~6]{BBL} implies that $f$ is differentiable (and thus pointwise Lipschitz)
at all $x\in\Rn\setminus A$ with $\Leb^n(A)=0$.
Lemma~\ref{ptwiselem} implies that $F$ is poitwise-Lipschitz at all $x\in\Rn\setminus A$.
\cite[Theorem~2.6]{Dmwd} implies that there exists $B\subset\Rn$ with $\Leb^n(B)=0$ such that
$F$ is metrically differentiable at all $x\in \Rn\setminus(A\cup B)$.
\par
If $M=Y^*$, then the result follows similarly from~\cite[Corollary~2.8]{Dmwd}.
\end{proof}

\begin{theorem}\label{kirthm}
Let $X$ be an separable Banach space, $K\subset X$ be a closed
convex cone with non-empty interior, and let $(M,\rho)$ be a metric space.
Then for every $K$-monotone dominated $F:X\to M$ there exists an Aronszajn
null set $A\subset X$ such that $F$ is metrically G\^ateaux differentiable at all points
of $X\setminus A$.
\end{theorem}

\begin{proof}
First, we can embed $(M,\rho)$ isometrically into
some $\ell_\infty(\Gamma)$, and so we can assume that $F:\Rn\to\ell_{\infty}(\Gamma)$.
Let $f:\Rn\to\R$ be a dominating function for $F$.
Then~\cite[Theorem~9]{BW} together with~\cite[Proposition~4]{BBL} implies that there exists 
an Aronszajn null set $A_1\subset X$ such that $f$ is G\^ateaux differentiable
at all $x\in X\setminus A_1$.
Lemma~\ref{pttoptlem} implies that $f$ is pointwise-Lipschitz at all $x\in X\setminus A_1$,
and thus
Lemma~\ref{ptwiselem} implies that $F$ is pointwise-Lipschitz at all $x\in X\setminus A_1$.
\cite[Theorem~5.4]{Dmwd} implies that there exists an Aronszajn null
set $A_2\subset X$  such that
$F$ is metrically G\^ateaux differentiable at all $x\in X\setminus(A_1\cup A_2)$.
\end{proof}

We have the following corollary, which follows from Theorem~\ref{kirthm}
in the same way as Corollary~\ref{domcor2} follows from
Theorem~\ref{domdifthm}.

\begin{corollary}\label{domcor1}
Let $X$ be an separable Banach space, $Y$ be an ordered normed linear
space such that $Y_+$ is well-based, $K\subset X$ be a closed
convex cone with non-empty interior.
Then for every $(K,Y_+)$-monotone mapping $F:X\to Y$ there exists an Aronszajn
null set $A\subset X$ such that $F$ is metrically G\^ateaux differentiable at all points
of $X\setminus A$.
\end{corollary}

Using~\cite[Theorem~5.3]{Dmwd} we can prove the following result
(we leave the details to the reader as the proof is similar
to the proof of Theorem~\ref{kirthm}).

\begin{theorem}\label{wskirthm}
Let $X$ be an separable Banach space, $K\subset X$ be a closed
convex cone with non-empty interior, and let $Y$ be a separable Banach space.
Then for every $K$-monotone dominated $F:X\to Y^*$ there exists an Aronszajn
null set $A\subset X$ such that for each $x\in X\setminus A$ we
have that $F$ is metrically G\^ateaux differentiable at~$x$,
$F$ is $w^*$-G\^ateaux differentiable at~$x$, and
$MD(F,x)(w)=\|wd(F,x)(w)\|$ for all $w\in X$.
\end{theorem}

\section{Acknowledgment} The author would like to thank Jan Rycht\'a\v{r}
for a useful discussion about renormings.

\end{document}